\documentclass[12pt]{article}
\author{Gun Srijuntongsiri\thanks{4163 Upson Hall, Cornell University, Ithaca, NY 14853. Email: gunsri@cs.cornell.edu.}, 
Stephen A. Vavasis\thanks{4130 Upson Hall, Cornell University, Ithaca, NY 14853. Email: vavasis@cs.cornell.edu.}}
\title{A Fully Sparse Implementation of a Primal-Dual Interior-Point Potential Reduction Method for 
Semidefinite Programming\thanks{Supported in part by NSF DMS 0434338 and NSF CCF 0085969.}}
\date{December 2, 2004}

\usepackage{amsmath, amsthm, amssymb,graphicx}

\begin{document}
\maketitle

\newtheorem{thm}{Theorem}[section]
\newtheorem{prop}{Proposition}[section]

\begin{abstract}
In this paper, we show a way to exploit sparsity in the problem data in 
a primal-dual potential reduction method for solving a class of semidefinite programs.  
When the problem data is sparse, the dual variable is also sparse, but the primal one 
is not. To avoid working with the dense primal variable, we apply 
Fukuda et al.'s theory of partial matrix completion and work with partial matrices instead.
The other place in the algorithm where sparsity should 
be exploited is in the computation of the search direction, where the gradient and 
the Hessian-matrix product of the primal and dual barrier functions must be computed in 
every iteration.  By using an idea from automatic differentiation in 
backward mode, both the gradient and the Hessian-matrix product can be computed in time 
proportional to the time needed to compute the barrier functions of sparse variables 
itself.  Moreover, the high space complexity that is normally associated with 
the use of automatic differentiation in backward mode can be avoided in this case.  
In addition, we suggest a technique to efficiently compute the determinant
of the positive definite matrix completion that is required to compute primal search 
directions. The method of obtaining one of the primal search directions that minimizes 
the number of the evaluations of the determinant of the positive definite completion is 
also proposed. We then implement the algorithm and test it on the problem of finding the 
maximum cut of a graph.
\end{abstract}

\section{Introduction}
Let $\mathcal{S}^n$ denote the space of $n \times n$ symmetric matrices. Given $A_p \in \mathcal{S}^n, p=1,2,\ldots,m, b \in \mathbb{R}^m,$ and $C \in \mathcal{S}^n$, semidefinite programming (SDP) problems in standard form are given as
\begin{equation}
\label{primal}
\begin{array}{ll}
\min_{X \in \mathcal{S}^n} & C \bullet X \\
$subject to:$ & A_p \bullet X = b_p, p = 1, \ldots, m, \\
 & X \geq 0.
\end{array}
\end{equation}
The notation $C \bullet X$ represents the inner product of $C$ and $X$, which is equal to $\sum_{ij} C_{ij}X_{ij}$.  The constraint $X \geq 0$ denotes that $X$ must be symmetric positive semidefinite, which means $w^T X w \geq 0$ for any $w \in 
\mathbb{R}^n$. Similar to linear programming, semidefinite programs have associated dual problems. The semidefinite program (\ref{primal}) is said to be in \emph{primal} form. Its \emph{dual}, which is also a SDP problem, is
\begin{equation}
\label{dual}
\begin{array}{ll}
\max_{y \in \mathbb{R}^m, S \in \mathcal{S}^n} & b^T y \\
$subject to:$ & \sum_{p=1}^m y_p A_p + S = C, \\ 
& S \geq 0,
\end{array}
\end{equation}
where $S$ is called the \emph{dual slack matrix}. 

Semidefinite programming has many applications in many fields (see \cite{boyd} for a list of applications). It can also be regarded as an extension of linear programming. As a result, various methods for solving linear programming have been extended to solve SDP. In particular, interior-point methods were first extended to SDP by Alizadeh \cite{alizadeh} and Nesterov and Nemirovskii \cite{nesterov} independently. The most effective interior-point methods are primal-dual approaches that use information from both primal and dual programs. Most primal-dual interior-point algorithms for SDP proposed fall into two categories: path-following and potential reduction methods. Our algorithm is of the potential reduction kind, in which a potential function is defined and each iterate reduces the potential by at least a constant amount. It is based on the primal-dual potential reduction method proposed by Nesterov and Nemirovskii in their book \cite{nesterov}.  

This paper focuses on the sparse case of SDP, where the data matrices $C$ and $A_p$'s consist of mostly zero entries. Because most problem data arising in practice are sparse, it is vital for an SDP solver to take advantage of the sparsity and avoid unnecessary computation on zero entries. The obstacle that prevents effective exploitation of sparsity in an SDP algorithm is that the primal matrix variable is dense regardless of the sparsity of the data.  To avoid this problem, Benson et al.\ proposed a pure dual interior-point method for sparse case \cite{benson}. Later, Fukuda et al.\ proposed a primal-dual algorithm using partial matrix and matrix completion theory to avoid the dense primal matrix \cite{fukuda}. Our algorithms follow Fukuda et al.'s suggestion and uses partial primal matrix to take advantage of sparsity in the primal-dual framework. In contrast, a recent work by Burer is also built upon Fukuda et al.'s idea of using partial matrix but his algorithm is a primal-dual path-following method based on a new search direction \cite{burer}.

In Nesterov and Nemirovskii's primal-dual potential reduction method, the computation of the search directions requires the gradient and Hessian-matrix product of the barrier functions. The currently common way to compute this gradient is not efficient in some sparse case. Our algorithm applies the idea from automatic differentiation in reverse mode to compute gradient and Hessian-matrix product in a more efficient manner for the sparse cases. Additionally, we suggest a technique that evaluates the barrier function value of a partial matrix efficiently in certain cases and an alternative way to compute the search directions when such evaluation is expensive. When the data matrices' aggregated
sparsity pattern forms a planar graph, our algorithm manages to 
reduce the time complexity to $O(n^{5/2})$ operations
and the space complexity to $O(n \log n)$ per SDP iterate. This is
a significant improvement from the $O(n^3)$ time complexity and the $O(n^2)$
space complexity per iterate of a typical SDP solver for planar case.

We start with review of necessary material on semidefinite programming, Nesterov and Nemirovskii's primal-dual potential reduction method, and sparse matrix computation in Section \ref{section_sdp}, \ref{section_nesterov}, and \ref{section_sparse}, respectively. Section \ref{section_dualnewton} covers the detail of our algorithm's computation of the dual Newton direction including how an idea from automatic differentiation in reverse mode is used to evaluate the gradient and Hessian-matrix product efficiently. The computation of the primal projected Newton direction as well as the efficient method in computing the determinant of the positive definite completion of a partial matrix are discussed in Section \ref{section_primalnewton}. Finally, the results of our algorithm on test instances of the problem of finding maximum cut are in Section \ref{section_exp}.

\section{Preliminaries on semidefinite programming}
\label{section_sdp}
We refer to $X$ and $(y,S)$ as \emph{feasible solutions} if they satisfy the constraints in (\ref{primal}) and (\ref{dual}), respectively. A \emph{strictly feasible solution} 
is a feasible solution such that $X$ (or $S$) are symmetric positive definite.  A matrix $A$ is 
\emph{symmetric positive definite} $(A > 0)$ if $w^T Aw > 0$ for any $w \in \mathbb{R}^n\setminus\{\mathbf{0}\}$.  
Problem (\ref{primal}) (resp. (\ref{dual})) is \emph{strictly feasible} if it contains a strictly feasible solution.

Let $\mathcal{P}$ (resp. $\mathcal{D}$) denote the set of feasible solutions of (\ref{primal}) (resp. (\ref{dual})), and $\mathcal{P'}$ (resp. $\mathcal{D'}$) denote the set of strictly feasible solutions of (\ref{primal}) (resp. (\ref{dual})). The \emph{duality gap}, which is the difference between primal and dual objective functions
\begin{eqnarray*}
C \bullet X - b^T y & = & \left(\sum_{i=1}^n y_i A_p + S \right) \bullet X - b^T y \\
& = & S \bullet X, \\
\end{eqnarray*}
is nonnegative at any feasible solution \cite{nesterov}. Under the assumption that (\ref{primal}) and (\ref{dual}) are strictly feasible and bounded, $(X^*,y^*,S^*)$ solves (\ref{primal}) and (\ref{dual}) if and only if $X^* \in \mathcal{P}, (y^*,S^*) \in \mathcal{D}$, and $S^* \bullet X^* = 0$.

\section{Nesterov and Nemirovskii's Primal-dual Potential Reduction Method}
\label{section_nesterov}
\emph{Primal-dual potential reduction methods} solve SDP programs by minimizing the potential function
\begin{equation}
\label{potentialfun}
\phi(X,S) = (n + \gamma \sqrt{n}) \ln (S \bullet X) - \ln \det X - \ln \det S : \mathcal{P} \times \mathcal{D} 
\rightarrow \mathbb{R},
\end{equation} 
where $\gamma > 0$ is a given constant parameter of the algorithm. Any sequence of iterates in which the potential $\phi(X,S)$ tends to negative infinity converges to (or at least has an accumulation point at) a strictly feasible solution $(X^*,y^*,S^*)$ such that $S^* \bullet X^* = 0$ and hence is optimal \cite{renegar}.  

Most primal-dual potential reduction methods begin at a strictly feasible iterate and compute the next strictly
feasible iterate while guaranteeing at least constant decrease in $\phi$ each iteration. The process is continued until an iterate with duality gap less than or equal to $\epsilon$ is found, where $\epsilon > 0$ is a given tolerance. Nesterov and Nemirovskii
proposed a polynomial-time primal-dual potential reduction method for more general convex programming problems in their 
book \cite{nesterov}, which is the basis of our algorithm. We now proceed to explain the ``large step" version of their
method as applied to the SDP (\ref{primal}) and (\ref{dual}). We call this method the ``Decoupled Primal-Dual" algorithm (DPD) since the computation of primal and dual directions are less directly coupled than in other primal-dual interior point methods. Given a current strictly feasible iterate $(X,S)$, compute
the next iterate $(X',S')$ as follows:

\begin{enumerate}

\renewcommand{\labelenumi}{(\roman{enumi})} 

\item Let 

\begin{displaymath}
M = \frac{(n + \gamma \sqrt{n})}{S \bullet X} S.
\end{displaymath}

Define the function

\begin{displaymath}
v(W) = -\ln \det W + M \bullet (W-X),
\end{displaymath}
where $W \in \mathcal{S}^n$.

\item Find the projected Newton direction $N$ of $v$ onto $\mathcal{P}'$ at X,

\begin{displaymath}
N = \textrm{argmin}_H\{v'(X) \bullet H + \frac{1}{2} (v''(X) H) \bullet H  :  A_p \bullet H = 0, p=1,2,\ldots,m\}.
\end{displaymath}

\item Let $\lambda$ be
\begin{displaymath}
\lambda = [(v''(X) N) \bullet N]^{1/2}.
\end{displaymath}

\item We then have

\begin{displaymath}
\Delta X_1 = \frac{N}{1+\lambda}, 
\end{displaymath}

and

\begin{displaymath}
\Delta S_1 = \frac{S \bullet X}{(n + \gamma \sqrt{n})}[-\nabla \ln\det(X)-(\nabla^2\ln\det(X))N ] - S
\end{displaymath}

as a primal and dual directions respectively.

\end{enumerate}

By swapping the roles of primal and dual in (i)-(iv), another pair of directions can be achieved as follow:

\begin{enumerate}
\addtocounter{enumi}{4} 
\renewcommand{\labelenumi}{(\roman{enumi})} 

\item Let 

\begin{displaymath}
\tilde{M} = \frac{(n + \gamma \sqrt{n})}{S \bullet X} X.
\end{displaymath}

Define the function

\begin{displaymath}
\tilde{v}(W) = -\ln \det W + \tilde{M} \bullet (W-S).
\end{displaymath}

\item Find the projected Newton direction $\tilde{N}$ of $\tilde{v}$ onto $\mathcal{D}'$ at $S$,
\begin{displaymath}
\tilde{N} = \textrm{argmin}_H\{\tilde{v}'(S) \bullet H + \frac{1}{2} (\tilde{v}''(S) H) \bullet H  :  
H = \sum_{p=1}^m z_p A_p, \textrm{ for some } z \in \mathbb{R}^m \}.
\end{displaymath}

\item Let $\tilde{\lambda}$ be
\begin{displaymath}
\tilde{\lambda} = [(\tilde{v}''(S) \tilde{N}) \bullet \tilde{N}]^{1/2}.
\end{displaymath}

\item Then

\begin{displaymath}
\Delta S_2 = \frac{\tilde{N}}{1+\tilde{\lambda}}, 
\end{displaymath}

and

\begin{displaymath}
\Delta X_2 = \frac{S \bullet X}{(n + \gamma \sqrt{n})}[-\nabla \ln\det(S)-(\nabla^2\ln\det(S))\tilde{N} ] - X
\end{displaymath}

are another dual and primal directions, respectively.

\item Find

\begin{displaymath}
\begin{array}{lll}
(h_1^*, h_2^*, k_1^*, k_2^*) & = & \textrm{argmin}_{h_1,h_2,k_1,k_2} \phi(X+h_1 \Delta X_1 + h_2 \Delta X_2, 
                                   S+k_1 \Delta S_1 + k_2 \Delta S_2) \\
$subject to:$  & & X+h_1 \Delta X_1 + h_2 \Delta X_2 > 0, $ and $ \\
& & S+k_1 \Delta S_1 + k_2 \Delta S_2 > 0. 
\end{array}
\end{displaymath}

\item Finally, set

\begin{eqnarray}
X' & = & X + h_1^* \Delta X_1 + h_2^* \Delta X_2, \nonumber \\
S' & = & S + k_1^* \Delta S_1 + k_2^* \Delta S_2. \nonumber
\end{eqnarray}

\end{enumerate}

Nesterov and Nemirovskii also showed that DPD algorithm achieves a constant reduction 
in $\phi$ at each iteration even when only two directions $\Delta X_1$ and 
$\Delta S_1$ are considered and potential minimization in step (ix) is not performed (that is, when fixing $h^*_1=k^*_1= 1$ and $h^*_2=k^*_2 =0$).

\section{Sparse matrix computation}
\label{section_sparse}
A \emph{sparse} matrix is a matrix with few nonzero entries.  Many problems' data encountered 
in practice are sparse.  By exploiting their structures, time and space required to perform operations
on them can be greatly reduced.  Many important applications of SDP, such as the problem of finding maximum cut, 
usually have sparse data, too.  For this reason, we consider the sparse case in this paper.  

To be able to discuss the exploitation of sparse data in SDP, background on chordal graph theory is needed, which is
addressed in the following section.

\subsection{Chordal graphs}

Let $G=(V,E)$ be a simple undirected graph.  A \emph{clique} of $G$ is a complete induced subgraph of $G$. A clique $C = (V',E')$ is \emph{maximal} if its vertex set $V'$ is not a proper subset of another clique. Let $Adj(v)=\{u \in V  :  \{u,v\} \in E\}$ denote the set of all vertices adjacent to a vertex $v \in V$.  A vertex $v$ is called \emph{simplicial} if all of its adjacent vertices $Adj(v)$ induce a clique.

For any cycle of $G$, a \emph{chord} is an edge joining two non-consecutive
vertices of the cycle.  Graph $G$ is said to be \emph{chordal} if each of its cycles of length 4 or greater has a chord. One 
fundamental property of a chordal graph is that it has a simplicial vertex, say $v_1$ and that the subgraph induced
by $V\setminus\{v_1\}$ is again chordal, which therefore has a simplicial vertex, say $v_2$.  By repeating this process, 
we can construct a \emph{perfect elimination ordering} of the vertices, say $(v_1,v_2,\ldots,v_n)$, such that 
$Adj(v_i) \cap \{v_{i+1},v_{i+2},\ldots,v_n \}$ induces a clique for each $i=1,2,\ldots,n-1$. It was shown by Fulkerson and Gross that a graph is chordal if and only if it has a perfect elimination ordering \cite{fulkerson}.

Given a perfect elimination ordering $(v_1,v_2,\ldots,v_n)$ of a chordal graph, its maximal cliques can be enumerated 
easily. A maximal clique containing the simplicial vertex $v_1$ is given by ${v_1} \cup Adj(v_1)$ and is unique. A
maximal clique not containing $v_1$ is a maximal clique of the chordal subgraph induced on $V\setminus \{v_1\}$. Therefore, by repeating this reasoning, the maximal cliques $C_r \subseteq V, r=1,2,\ldots,l$, are given by 
\[
C_r = \{v_i\} \cup (Adj(v_i) \cap \{ v_{i+1},v_{i+2},\ldots,v_n \})
\] for $i = \min \{j  :  v_j \in C_r\}$, that is, the maximal members of $\{ \{v_i\} \cup (Adj(v_i) \cap \{ v_{i+1},v_{i+2},\ldots,v_n \}) : i=1,2,\ldots,n \}$.

One property of the sequence of maximal cliques is that it can be reindexed such that for any $C_r, r=1,2,\ldots,l-1$, there exists a $C_s, s \geq r+1$, such that
\begin{equation}
\label{running}
C_r \cap (C_{r+1} \cup C_{r+2} \cup \cdots \cup C_l ) \subsetneq C_s.
\end{equation} Such property is called the \emph{running intersection property}.

There is a well-known relationships between chordal graph and Cholesky factorization of sparse symmetric positive definite matrices. Given a symmetric positive definite matrix $X$, its Cholesky factor $L$ is a lower-triangular matrix such that $X=LL^T$. 
The \emph{sparsity pattern} of $X$, which is defined as the set of row/column indices of nonzero entries of $X$, is often represented as a graph $G = (V,E)$, where $V = \{1,2,\ldots,n\}$ and $E = \{\{i,j\} : 
X_{ij} \neq 0, i \neq j \}$.  Similarly, the sparsity pattern of $L$ can be represented by the graph $G'=(V,F)$, where
$F=\{\{i,j\} : L_{ij} \neq 0, i \geq j \}$. Under \emph{no numerical cancellations assumption}, which means no zero entries are resulted from arithmetic operations on nonzero values, it is seen that $F \supseteq E$, with $F$ having possibly additional fill-ins. In addition, $G'=(V,F)$ is chordal and is said to be a \emph{chordal extension} of $G=(V,E)$.  

The number of fill-ins in the Cholesky factorization depends on the ordering of the row/column indices. The question of
finding the reordering of the row/column indices to yield fewest fill-ins is NP-complete.  In the best case when $G$ is chordal, a perfect elimination ordering yields the Cholesky factor with no fill-ins.

\subsection{Partial symmetric matrix and positive definite matrix completion}

It is a well-known fact that in the course of SDP algorithms, the primal variable $X$ usually is dense 
even if the data are sparse while the dual variables $y$ and $S$ stay sparse.  
To avoid working with a dense primal variable, Fukuda et al.\ suggested the use of a partial symmetric matrix
for the primal variable in SDP algorithms \cite{fukuda}. Let $V = \{1,2,\ldots,n\}$. Define the 
\emph{aggregate sparsity pattern $E$} of the data to be
\begin{displaymath}
E = \{(i,j) \in V \times V  :  C_{ij} \neq 0 \textrm{ or } [A_p]_{ij} \neq 0 \textrm{ for some } p \in 
\{0,1,\ldots,m\}\}.
\end{displaymath}

Observing (\ref{primal}), we see that the values of the objective function and constraint linear functions
only depend the entries of $X$ corresponding to the nonzero entries of $C$ and $A_p$'s.  The remaining entries 
of $X$ affect only whether $X$ is positive semidefinite.  In other words, if $X$ and $X'$ satisfy $X_{ij} = X'_{ij}$,
for any $(i,j) \in E$, then

\begin{eqnarray}
C \bullet X & = & C \bullet X' \nonumber \\
A_p \bullet X & = & A_p \bullet X', p = 1, 2, \ldots, m. \nonumber
\end{eqnarray}

A \emph{partial symmetric matrix} is a symmetric matrix in which not all of its entries are specified.  
A partial symmetric matrix $\bar{X}$ can be treated as a sparse matrix, having its unspecified entries regarded as 
having zero values. Hence, a sparsity graph $G'=(V,F)$ can be used to represent the row/column indices of specified
entries of $\bar{X}$ in the same manner as it is used to represent nonzero entries of a sparse matrix.  Let 
$\mathcal{S}^n(F,?)$ denote the set of $n \times n$ partial symmetric matrices with entries specified in $F$. We assume that all diagonal entries are also specified although there are no edges in $G'$ representing them.

A \emph{completion} of a partial symmetric matrix $\bar{X}$ is a matrix $X$ of the same size as $\bar{X}$ such that
$X_{ij} = \bar{X}_{ij}$ for any $\{i,j\} \in F$.  A \emph{positive definite completion} of a partial symmetric matrix is a completion that is positive definite. The following theorem characterizes when a partial matrix has a positive definite completion. 

\begin{thm}[{Grone et al.\ \cite[Theorem 7]{grone}}]
\label{completion_theorem} 
Let $G'=(V,F)$ be a chordal graph. Any partial symmetric matrix $\bar{X} \in \mathcal{S}^n(F,?)$ satisfying the property that $\bar{X}_{C_r C_r}$ is symmetric positive definite for each $r=1,2,\ldots,l$, where $\{ C_r \subseteq V : r = 1,2,\ldots,l\}$ denote the family of maximal cliques of $G'$, can be completed to a positive definite matrix.
\end{thm}

\subsection{Maximum-determinant positive definite matrix completion}
\label{section_maxdetcompletion}

The following result of Fukuda et al.\ \cite{fukuda} shows an efficient way to compute a certain positive definite matrix completion. Given a partial symmetric matrix $\bar{X}$ whose sparsity pattern $G'=(V,F)$ is chordal, its unique positive definite completion that maximizes the determinant 
\[
\hat{X} = \textrm{argmax}_X \{\det (X) : X \textrm{ is a positive definite completion of } \bar{X} \}
\] is shown to be 
\begin{equation}
\label{cliquefactor}
P \hat{X} P^T = L_1^T L_2^T \cdots L_{l-1}^T D L_{l-1} \cdots L_2 L_1,
\end{equation} where $P$ is the permutation matrix such that $(1,2,\ldots,n)$ is the perfect elimination ordering for 
$P \bar{X} P^T$, $L_r$ ($r=1,2,\ldots,l-1)$ are sparse triangular matrices, and $D$ is a positive definite block-diagonal
matrix, both defined below \cite{fukuda}. Let $(C_1,C_2,\ldots,C_l)$ be an ordering of maximal cliques of $G'$ that enjoys the running intersection property (\ref{running}). Define
\[
\begin{array}{llll}
S_r & = & C_r \setminus (C_{r+1} \cup C_{r+2} \cup \cdots \cup C_l), &  r = 1,2,\ldots,l, \\
U_r & = & C_r \cap (C_{r+1} \cup C_{r+2} \cup \cdots \cup C_l), & r = 1,2,\ldots,l.
\end{array}
\] The factors in (\ref{cliquefactor}) are given by
\begin{equation}
[L_r]_{ij} = \left\{\begin{array}{ll}
         1, & i=j,\\
        {[\bar{X}^{-1}_{U_r U_r} \bar{X}_{U_r S_r}]}_{ij}, & i \in U_r, j \in S_r, \\
         0, & \textrm{otherwise}\\
      \end{array} \right.
\end{equation} for $r=1,2,\ldots,l-1$, and
\[
D= \left( \begin{array}{cccc}
    D_{S_1 S_1} & & & \\
    & D_{S_2 S_2} & & \\
    & & \ddots & \\
    & & & D_{S_l S_l} \\
    \end{array} \right),
\] where
\begin{equation}
\label{dsrsr}
D_{S_r S_r} = \left\{ \begin{array}{ll}
              \bar{X}_{S_r S_r} - \bar{X}_{S_r U_r}\bar{X}^{-1}_{U_r U_r} \bar{X}_{U_r S_r}, & r = 1,2,\ldots,l-1, \\
              \bar{X}_{S_l S_l}, & r = l.
              \end{array} \right.
\end{equation}
In addition, the unique determinant-maximizing positive definite completion $\hat{X}$ has the property that 
\begin{equation}
\label{xinvprop}
\begin{array}{rlll}
(\hat{X}^{-1})_{ij} &= &0, & (i,j) \notin F.
\end{array}
\end{equation}
In other words, the inverse of the determinant-maximizing completion has the same sparsity pattern as that of the partial matrix.

\subsection{Using the maximum-determinant extension of $\bar{X}$ in DPD}

To exploit sparsity in the data matrices, our algorithm works in the space of 
partial matrix $\bar{X}$ for primal variable. When positive definite completion 
of $\bar{X}$ is needed, the maximum-determinant completion is used. 
We choose this particular completion because it preserves the self-concordance 
of the barrier function, which follows directly from Proposition 5.1.5 of \cite{nesterov}.
This property guarantees that our algorithm converges to an optimal solution.

\section{Computation of the dual Newton direction}
\label{section_dualnewton}

The single most computationally-intensive step in any potential reduction method is the computation of search directions. In Nesterov and Nemirovskii's method described in Section \ref{section_nesterov}, this computation occurs in steps (ii), (iv), 
(vi), and (viii). For this reason, minimizing computation in these steps are emphasized in our algorithm.  We describe our algorithm to compute $\tilde{N}$ (step (vi) of the algorithm), which is the most computationally-intensive step in the computation of $\Delta S_2$, in this section. Conjugate gradient is used together with an idea from automatic differentiation in reverse mode to compute $\tilde{N}$.

The minimization problem in step (vi) is
\begin{equation}
\label{minNt}
\begin{array}{ll}
\min_{\tilde{N}} & -(\nabla\ln\det S+\tilde{M}) \bullet \tilde{N} - \frac{1}{2}((\nabla^2\ln\det S)\tilde{N}) \bullet \tilde{N} \\
$subject to:$ & \tilde{N} = \sum_{p=1}^m z_pA_p, \textrm{ for some } z \in \mathbb{R}^m.
\end{array}
\end{equation}
Replacing $\tilde{N}$ with $\sum_{p=1}^m z_pA_p$ in (\ref{minNt}) yields 
\begin{displaymath}
\min_{z} -(\nabla\ln\det S+\tilde{M}) \bullet \sum_{p=1}^m z_pA_p - \frac{1}{2}((\nabla^2\ln\det S)\sum_{p=1}^m z_pA_p) \bullet \sum_{p=1}^m z_pA_p,
\end{displaymath}
which is equivalent to solving the system
\begin{equation}
\label{eqz}
\mathcal{A}(\sum_{p=1}^m z_p(\nabla^2\ln\det S) A_p) = -\mathcal{A}(\nabla\ln\det S-\tilde{M})
\end{equation}
for $\mathbf{z}$, where $\mathcal{A}(W)= \left(\begin{array}{c}A_1\bullet W\\ \vdots\\A_m\bullet W\end{array}\right)$. To rewrite (\ref{eqz}) into standard system of linear equations form, we first define the functions
\begin{displaymath}
f(\mathbf{u}) \equiv \ln\det(C-\sum_{p=1}^mu_pA_p)
\end{displaymath}
and
\[
h(\mathbf{u}) \equiv f(\mathbf{u}) + \tilde{M} \bullet (C-\sum_{p=1}^mu_pA_p-S).
\]
Then, the system (\ref{eqz}) is equivalent to 
\begin{equation}
\label{eqz2}
(\nabla^2h(\mathbf{y}))\mathbf{z} = -\nabla h(\mathbf{y}).
\end{equation}
Our algorithm uses conjugate gradient to solve (\ref{eqz2}) to exploit the fact that $(\nabla^2h(\mathbf{y}))\mathbf{z}$ and $\nabla h(\mathbf{y})$ can be computed efficiently.

\subsection{Computing derivatives}
\label{section_deriv}

As described above, the conjugate gradient method when solving for $\tilde{N}$ 
calls for the computation of $(\nabla^2h(\mathbf{y}))\mathbf{z}$ in each 
iteration and $\nabla h(\mathbf{y})$ once. Note first that $h(\mathbf{u})$ can be 
can be computed as follow: Cholesky factorize $C-\sum_{p=1}^m u_p A_p = LL^T$, 
where $L$ is a lower triangular matrix, 
compute $\ln\det(C-\sum_{p=1}^m u_p A_p) = \ln\det(LL^T) = 2\ln\det L= 2\sum_i\ln L_i$,
and finally $h(\mathbf{u}) = 2\sum_i\ln L_i + \tilde{M} \bullet (C-\sum_{p=1}^m u_p A_p-S)$.
We derive the algorithm to evaluate $\nabla h(\mathbf{u})$ from the above method of
evaluating $h(\mathbf{u})$ by imitating automatic differentiation (AD) in reverse mode,
which is discussed in detail below. 
To evaluate $(\nabla^2h(\mathbf{u}))\mathbf{z}$, 
notice that it is the derivative of the function $g(\mathbf{u}) \equiv [h(\mathbf{u})]^T\mathbf{z}$ with respect to $\mathbf{u}$. Hence, we can derive the algorithm to evaluate 
$(\nabla^2h(\mathbf{u}))\mathbf{z}$ from the algorithm to compute $g(\mathbf{u})$, again by
imitating AD in reverse mode. We emphasize that we do not suggest using AD to automatically
compute derivatives of $h(\mathbf{u})$ given the algorithm to evaluate $h(\mathbf{u})$ as
we would not be able to control the space allocation of AD. Rather, we imitate how 
AD in reverse mode differentiate the algorithm to evaluate $h(\mathbf{u})$, make
additional changes to reduce space requirement (discussed below), and then 
hand-code the resulting algorithm.

\emph{Automatic differentiation} is a tool that receives a code that evaluates a function as its input and generates a new piece of code that computes the value of the first derivative of the same function at a given point in addition to evaluating the function.  In essence, AD repeatedly applies the chain rule to the given code.  There are two modes in AD, each representing a different approach in applying the chain rule. \emph{Forward mode} differentiates each intermediate variable with respect to each input variable from top to bottom. \emph{Reverse mode}, on the other hand, differentiates each output variable with respect to each intermediate variable from bottom up, hence the name \emph{reverse}.  Note that each entry in a matrix is treated individually. Therefore, one $n \times n$ input matrix is treated as $n^2$ input variables \cite{griewank}.

One mode is more suitable than the other in different situations. Complexity-wise, forward mode is more appealing when the number of input variables is less than the number of output variables while reverse mode is more appealing when the number of input variables is greater. Let $\omega(f)$ be the computation time of the given code, $c$ be the number of input variables of the code, and $d$ be the number of output variables.  The code generated by forward mode computes the first derivative in time proportional to $c \omega(f)$ while the one generated by reverse mode does so in time proportional to $d \omega(f)$. However, reverse mode has one additional disadvantage: the storage space required may be as large as time complexity of the original code, which can be much larger than the space complexity of the original code, because if a variable is updated many times throughout the evaluation of $f$, its values before and after each such update may be needed.  Forward mode does not suffer from this problem because by taking derivatives from top to bottom, the old value of an intermediate variable is not needed after the variable is updated and therefore can be safely overwritten in the same storage space.  The storage issue in reverse mode can be partially fixed by recomputing required values rather than storing them, but this approach may result in significant increase in computation time.

For our problem, however, reverse mode can be applied to compute 
$\nabla h(\mathbf{y})$ and $(\nabla^2h(\mathbf{y}))\mathbf{z}$
without increasing storage requirement.  By performing reverse mode AD by hand,
it is seen that all of the intermediate variables can be overwritten safely 
and thus avoiding the need to store many versions of a variable.  
Therefore, our method of computing $\nabla h(\mathbf{y})$ and 
$(\nabla^2h(\mathbf{y}))\mathbf{z}$ requires the same order of time and space
complexity as the algorithm for evaluating the original function $h(\mathbf{y})$.

Analytically, it can be shown that $\nabla f(\mathbf{y}) = \mathcal{A}(S^{-1})$.
From the definition of $\mathcal{A}(\cdot)$, we see that only entries of 
$S^{-1}$ in $F$, the chordal extension of the aggregated sparsity pattern $E$, 
need to be computed in order to compute $\nabla f(\mathbf{y})$ (and, consequently,
$\nabla h(\mathbf{y})$). Erisman and Tinney showed a method of computing such
entries of $S^{-1}$ in the same order of time and space complexity as performing 
Cholesky factorization of $S$ in 1975 \cite{erisman}. Thus, their method can
be used to compute $\nabla h(\mathbf{y})$ in the same complexity as our proposed
method. Nevertheless, our method proves useful as it can be extended to compute 
the Hessian-vector product efficiently.

This idea of imitating reverse AD is not limited to computing derivatives of 
$h(\mathbf{u})$. It can also be applied to compute the gradient of $\ln\det S$ 
with respect to entries in $F$ of $S$, which is required in step (vi) and (viii)
of our algorithm. The most computationally expensive part of evaluating 
$\ln\det S$ is to Cholesky factorize $S$, which is similar to the algorithm 
for evaluating $h(\mathbf{u})$. Hence, their derivative codes are very similar,
and all of the intermediate variables arising from performing reverse mode AD
on $\ln\det S$ evaluation algorithm can be safely overwritten, too. 
The entries in $F$ of $(\nabla^2\ln\det S)\tilde{N}$ required in step (viii) 
can also be computed using the same idea since $(\nabla^2\ln\det S)\tilde{N} =
\frac{d}{dS}\left(\left(\nabla \ln\det S\right)\bullet \tilde{N}\right)$.

We remark that our approach can be much more efficient than 
the obvious way of obtaining $\nabla h(\mathbf{u})$ or the gradient of 
$\ln\det S$ with respect to entries in $F$. The simple way of obtaining 
the entries in $F$ of $S^{-1}$ is to (i) compute the Cholesky factorization 
$S = L L^T$ and then (ii) compute the required entries of $S^{-1}$ by using 
backward and forward substitution to solve linear systems of the form 
$LL^T v_i = e_i$ for the required entries of $v_i$, the $i$th column of $S^{-1}$, 
where $e_i$ is the $i$th column of the identity matrix.  
When $S$ is sparse, step (ii) may be much more computationally expensive than step (i).  
One example is when $S$ is tridiagonal, in which case, step (i) requires only $O(n)$ 
operations while step (ii) requires $O(n)$ operations per entry of $S^{-1}$, 
which can result in a total of $O(n^2)$ operations if the number of nonzeros in $W$ 
is $O(n)$. On the other hand, our algorithm would require only $O(n)$ operations in this
case.

Although it appears that the sparse-inverse algorithm has not been previously used
in semidefinite programming, it has been used elsewhere in the optimization
literature. See, for example, Neumaier and Groeneveld \cite{neumaier}.

\section{Computation of the primal projected Newton direction}
\label{section_primalnewton}

Following the discussion in Section \ref{section_sparse}, our algorithm works with a partial matrix $\bar{X}$ with specified entries in $F$, a chordal extension of the aggregated sparsity pattern $E$, for primal variable.  For this reason, the computation of the primal search directions are more complicated than the dual ones described in previous section. Moreover, as we shall see below, the evaluation of $(\nabla^2 \ln\det \hat{X})P$, where $\hat{X}$ is the maximum determinant positive definite completion of $\bar{X}$ and $P$ is an arbitrary matrix, appears to be more expensive than performing Cholesky factorization. Consequently, the same algorithm used to compute the dual Newton direction as described in Section \ref{section_dualnewton}, which involves evaluation of $(\nabla^2 \ln\det \hat{X})P$ in each iteration of the conjugate gradient, may not be efficient. Therefore in this section, we propose a different method for obtaining $N$ that avoids excessive evaluation of $(\nabla^2 \ln\det \hat{X})P$.

Assume $\hat{X}^{-1}$ is known in addition to $\bar{X}$ (the detail on the computation of $\hat{X}^{-1}$ is addressed below in Section \ref{section_xderiv}). Recall from (\ref{xinvprop}) that $\hat{X}^{-1}$ has sparsity pattern $F$ and therefore is sparse. To compute for $N$, according to step (ii), the problem under consideration is
\begin{equation}
\label{minN}
\begin{array}{ll}
\min_{N} & -(\nabla\ln\det \hat{X}+M) \bullet N - \frac{1}{2}((\nabla^2\ln\det \hat{X})N) \bullet N \\
$subject to:$ & A_p \bullet N = 0, p = 1, 2, \ldots, m.
\end{array}
\end{equation}
Note that $\nabla\ln\det\hat{X}$ is $\hat{X}^{-1}$ and $(\nabla^2\ln\det \hat{X})N$ is $\hat{X}^{-1}N\hat{X}^{-1}$. The KKT condition for the optimum solution to (\ref{minN}) is
\begin{displaymath}
\hat{X}^{-1}N\hat{X}^{-1} = \hat{X}^{-1} - M + \sum_{p=1}^m \lambda_pA_p,
\end{displaymath}
or, equivalently, 
\begin{equation}
\label{lag1}
N = \hat{X} - \hat{X}M\hat{X} + \sum_{p=1}^m \lambda_p\hat{X}A_p\hat{X},
\end{equation}
where $\lambda_p$ $(p=1,2,\ldots,m)$ is a scalar to be determined that enforces the condition $A_p \bullet N = 0$ $(p=1,2,\ldots,m)$. To determine $\lambda_p$'s, eliminate $N$ from (\ref{lag1}) by taking inner product with $A_q$ $(q=1,2,\ldots,m)$ on both sides and noting that $A_q \bullet N = 0$, yielding the linear system 
\begin{equation}
\label{lagnoN}
\mathcal{A}(\sum_{p=1}^m \lambda_p\hat{X}A_p\hat{X}) = \mathcal{A}(\hat{X}M\hat{X} -\hat{X}),
\end{equation}
where $\mathcal{A}(\cdot)$ is defined as in Section \ref{section_dualnewton}. 
To rewrite (\ref{lagnoN}) as a standard system of linear equations form, define the function
\[
q(\mathbf{u}) = \ln\det(\hat{X}^{-1} - \sum_{p=1}^mu_pA_p).
\]
The system (\ref{lagnoN}) is therefore equivalent to
\begin{equation}
\label{lag2}
(\nabla^2 q(\mathbf{0}))\mathbf{\lambda} = \mathcal{A}(\hat{X}M\hat{X} -\hat{X}),
\end{equation}
where $\mathbf{\lambda} = (\lambda_1,\lambda_2,\ldots,\lambda_m)^T$. Conjugate gradient is then used to solve the system (\ref{lag2}) for $\mathbf{\lambda}$. After knowing $\mathbf{\lambda}$, we can now compute $N$ from (\ref{lag1}).

To solve for $\lambda_p$ efficiently with conjugate gradient, it is important that $(\nabla^2 q(\mathbf{0}))\mathbf{\lambda}$ and $\hat{X}M\hat{X}$ are not expensive to evaluate. This is where $\hat{X}^{-1}$ becomes useful. From the definition of $\mathcal{A(\cdot)}$, we see that we do not need to know the entries outside $F$ of the resulting matrices $\hat{X}M\hat{X}$. Also, the matrix $\hat{X}^{-1}$, unlike $\bar{X}$, is not a partial matrix, but recall from Section \ref{section_maxdetcompletion} that $\hat{X}^{-1}$ has the same sparsity pattern $F$ as the partial matrix $\bar{X}$. Moreover, $\hat{X}M\hat{X} = \left. \frac{d^2}{dW^2} \left((\ln\det W)M\right) \right\arrowvert_{W=\hat{X}^{-1}}$ (see appendix C of \cite{klerk}). Therefore, the entries in $F$ of $XMX$ can be computed using the idea of automatic differentiation in reverse mode in the same manner as computing $(\nabla^2\ln\det S)P$, as detailed in Section \ref{section_deriv}. The Hessian-vector product $(\nabla^2 q(\mathbf{0}))\mathbf{\lambda}$ can also be handled in the same manner as $(\nabla^2 h(\mathbf{y}))\mathbf{z}$ in the dual case. 
Lastly, the term $\hat{X}$ that is by itself in the quantity $(\hat{X} - \hat{X}M\hat{X})$ of (\ref{lag2}) may be replaced by $\bar{X}$ safely as the entries outside $F$ of $\hat{X}$ do not affect the equation after the inner product with $A_q$ is taken.

\subsection{Logarithm of determinant of positive definite completion matrix}
\label{section_maxdet}

Steps (i)-(ii) and (v)-(viii) of DPD can be performed using the techniques described 
in previous sections and the values of the matrices $S$, $\bar{X}$, and $\hat{X}^{-1}$. 
For step (ix), steepest descent method used to compute step size requires that the algorithm 
evaluates 
$\phi(X+h_1 \Delta X_1 + h_2 \Delta X_2,S+ k_1 \Delta S_1 + k_2 \Delta S_2)$ for a 
current point $(h_1,h_2,k_1,k_2)$ to be able to decide when to terminate the steepest descent.
But since we only have the partial matrices $\bar{X}$, $\overline{\Delta X}_1$, and $\overline{\Delta X}_2$,
we need to to be able to evaluate $\ln\det\hat{X}$ after we update $\bar{X}$ as 
$\bar{X}+h_1 \overline{\Delta X}_1 + h_2 \overline{\Delta X}_2)$.
Computation of $\ln\det \hat{X}$ is not trivial because, unlike 
the objective function or the linear constraints, the value of $\ln\det \hat{X}$ does depend on the 
entries outside the aggregated sparsity pattern $F$. 
We cover an efficient algorithm to compute $\ln\det \hat{X}$ in this section.

Consider a partial symmetric matrix $\bar{X}$ with sparsity pattern $F$, a chordal extension of $E$. Using the factors given in (\ref{cliquefactor}), the value of $\ln \det \hat{X}$ can be evaluated efficiently as follows. 
Because each $L_r (r=1,2,\ldots,l-1)$ is unit lower triangular, its determinant is one. The determinant of the block diagonal
matrix $D$ is the product of the determinants of each of its diagonal blocks $D_{S_r S_r}, r=1,2,\ldots,l$. Observe that 
$D_{S_r S_r} = \bar{X}_{S_r S_r} - \bar{X}_{S_r U_r}\bar{X}^{-1}_{U_r U_r} \bar{X}_{U_r S_r} (r=1,2,\ldots,l-1)$ in (\ref{dsrsr}) is the Schur complement of $\bar{X}_{U_r U_r}$ in
\[
Q \bar{X}_{C_r C_r} Q^T = \left( \begin{array}{cc}
                            \bar{X}_{S_r S_r} & \bar{X}_{S_r U_r} \\
                            \bar{X}_{U_r S_r} & \bar{X}_{U_r U_r}
                            \end{array} \right),
\] for some permutation matrix $Q$. The determinant of the Schur complement is
\begin{eqnarray}
\det (D_{S_r S_r}) & = & \frac{\det(Q \bar{X}_{C_r C_r} Q^T)}{\det(\bar{X}_{U_r U_r})} \nonumber \\
                   & = & \frac{\det(\bar{X}_{C_r C_r})}{\det(\bar{X}_{U_r U_r})}, \nonumber
\end{eqnarray} for $r = 1,2,\ldots,l-1$. Therefore,
\begin{eqnarray}
\ln\det \hat{X} & = & \ln\det (P\hat{X}P^T) \nonumber \\
                & = & \ln\left(\frac{\prod_{r=1}^l \det(\bar{X}_{C_r C_r})}{ \prod_{r=1}^{l-1} 
                      \det(\bar{X}_{U_r U_r})}\right) \nonumber \\
                & = & \sum_{r=1}^l \left(\ln\det \bar{X}_{C_r C_r}\right) - \sum_{r=1}^{l-1} \left(\ln\det \bar{X}_{U_r U_r}\right).
\label{logmaxdet}
\end{eqnarray}

Note that performing Cholesky factorization of a positive definite matrix with sparsity pattern $F$ requires $O(\sum_{i=1}^n s_i^2)$, where $s_i = \max_{C_r \ni i} |C_r \cap \{i,i+1,\ldots,n\}|$, assuming $(1,2,\ldots,n)$ is a perfect elimination ordering. On the other hand, computing $\ln\det\hat{X}$ by straightforward application of (\ref{logmaxdet}), that is, by computing determinants of each $\bar{X}_{C_r C_r}$ and $\bar{X}_{U_r U_r}$ separately, requires $O(\sum_{r=1}^l |C_r|^3)$ operations. Notice that the time required to perform Cholesky factorization of a positive definite matrix with sparsity pattern $F$ is the lower bound of the computation time of $\ln\det\hat{X}$, which occurs when $\hat{X}=\bar{X}$. For this reason, we seek to find an algorithm that computes $\ln\det\hat{X}$ in the same order of complexity as that of performing Cholesky factorization on the same sparsity pattern.

In the most favorable case where none of the maximal cliques overlap, straightforward application of (\ref{logmaxdet}) has the same time complexity as that of Cholesky factorization. To see the equivalence of the two algorithms' complexity in this case, note that $O(\sum_{i=1}^n s_i^2) = O(\sum_{r=1}^l \sum_{i \in C_r} s_i^2) = O(\sum_{r=1}^l |C_r|^3)$. An example of such case is when $\bar{X}$ is block diagonal. 

We consider the efficiency of straightforward calculation of (\ref{logmaxdet}) in the case that the sparsity pattern graph $G=(V,E)$ is planar next as this special case arises often in practice.  Our analysis assumes that the vertices of $G$ are ordered according to the nested dissection ordering. Lipton et al.\ introduce generalized nested dissection and show that performing Cholesky factorization on said ordering requires $O(n^{3/2})$ operations, where $n$ is the number of vertices \cite{lipton}. Planar graphs satisfy a \emph{$\sqrt{n}$-separator theorem}, which states that the vertices of the graph $G$ can be partitioned into three sets $A, B,$, and $C$ such that there are no edges having one endpoint in $A$ and the other in $B$, $|A|, |B| \leq \frac{2}{3}n$, and $|C| \leq \sqrt{8n}$. \emph{Nested dissection ordering} is computed by partitioning $V$ into $A$,$B$, and $C$ according to the separator theorem, number the unnumbered vertices in $C$ such that they are eliminated after the unnumbered vertices in $A$ and $B$, and then recursively number the unnumbered vertices in $A \cup C$ and $B \cup C$. The recursion stops when the number of vertices under consideration is less than $72$, at which point, the unnumbered vertices are numbered arbitrarily.  It is shown in Lipton et al.\ that, for a given $A, B,$ and $C$ in any level of the recursion, no vertex in $A$ is adjacent to any vertex in $B$ in the chordal extension graph $G'$. Therefore, any maximal clique of $G'$ can contain at most the vertices in the separator $C$ of each recursion hierarchy and additional $72$ vertices from the lowest level of recursion. Since each recursion reduces the number of vertices to at most $\frac{2}{3}n'$, where $n'$ is the number of vertices in consideration of the current level, and the separator has at most $\sqrt{8n'}$ vertices, the number of vertices in any maximal clique is at most $\sqrt{8n} + \sqrt{8(\frac{2}{3}n)} + \sqrt{8((\frac{2}{3})^2n)}+\ldots + \sqrt{8((\frac{2}{3})^{\log_{3/2}\frac{n}{72}}n)} + 72 = \sqrt{8n}\left(1+\sqrt{\frac{2}{3}} + \left(\sqrt{\frac{2}{3}}\right)^2 + \ldots + \left(\sqrt{\frac{2}{3}}\right)^{\log_{3/2}\frac{n}{72}}\right) + 72 = O(\sqrt{n})$. The number of maximal cliques is no greater than $n$. Therefore, straightforward calculation of $\ln\det\hat{X}$ requires $O(n (\sqrt{n})^3) = O(n^{5/2})$ operations, which is greater than that of Cholesky factorization by a factor of $n$. 

As seen from the planar case, straightforward computation of $\ln\det\hat{X}$ can be significantly more expensive than Cholesky factorization. Another common case that suffers from the same problem is when $\bar{X}$ is a banded matrix. A \emph{banded matrix} with \emph{bandwidth $p$} satisfies the property that the entry $\bar{X}_{ij} = 0$ if $|i-j| > p$. Performing Cholesky factorization on such a matrix takes $O(np^2)$ operations. To analyze time complexity of $\ln\det\hat{X}$ computation, first notice that $(1,2,\ldots,n)$ is a perfect elimination ordering and that the sequence of maximal cliques $\{ C_1, C_2,\ldots,C_{n-p} \}$, where $C_r = \{r,r+1,\ldots,r+p \}$ $(r=1,2,\ldots,n-p)$, satisfies the running intersection property (\ref{running}). Therefore, straightforward computation of $\ln\det\hat{X}$ requires $O(\sum_{r=1}^l |C_r|^3) = O((n-p)(p+1)^3) = O(np^3)$ operations, which is greater than $O(np^2)$ operations of Cholesky factorization.

However, it is possible to reduce the complexity of computing $\ln\det\hat{X}$ to $O(np^2)$ operations in the banded matrix case by using the following idea. The determinant of each (positive definite) submatrix $\bar{X}_{C_r C_r}$ and $\bar{X}_{U_r U_r}$ is usually computed from the product of diagonal entries of its Cholesky factor.  If a set $C_r$ (resp. $U_r$) shares many members with another set $C_{r'}$ (resp. $U_{r'}$), $r \neq r'$, the Cholesky factor of a symmetric permutation of $\bar{X}_{C_{r'} C_{r'}}$ (resp. $\bar{X}_{U_{r'} U_{r'}}$) can be constructed from the Cholesky factor of $\bar{X}_{C_r C_r}$ (resp. $\bar{X}_{U_r U_r}$) or \emph{vice versa}, which is more efficient than computing Cholesky factor of $\bar{X}_{C_{r'} C_{r'}}$ (resp. $\bar{X}_{U_{r'} U_{r'}}$) from scratch. In the banded matrix case, any two adjacent cliques $C_r$ and $C_{r+1}$ $(r=1,2,\ldots,n-p-1)$ share the same $p-1$ elements $(C_r \cap C_{r+1} = \{r+1,r+2,\ldots,r+p\})$. The same can be said about adjacent $U_r$'s. Observe that $U_r = \{r+1,r+2,\ldots,r+p\}$ $(r=1,2,\ldots,n-p-1)$ and $U_{n-p} = \emptyset$. Therefore, the two adjacent $U_r$ and $U_{r+1}$ $(r=1,2,\ldots,n-p-2)$ share $p-2$ elements.

The process of updating a Cholesky factor is as follow. Let $L_r$ be the Cholesky factor of $\bar{X}_{C_r C_r}$. Remove the first row of $L_r$, which corresponds to the $r$th row/column of $\bar{X}$, and let $\tilde{L}_r$ be the resulting $(p-1) \times p$ submatrix. We then transform $\tilde{L}_r^T$ to a $(p-1) \times (p-1)$ upper triangular matrix $R$ by performing Givens rotations to zero out the $p-1$ entries below the main diagonal of $\tilde{L}_r^T$. Notice that the columns of $R$ corresponds to the $(r+1)$th, $(r+2)$th,\ldots,$(r+p-1)$th columns of $\bar{X}$. Therefore, $R^T$ is exactly the first $p-1$ rows of the Cholesky factor $L_{r+1}$ of $\bar{X}_{C_{r+1} C_{r+1}}$. The final row of $L_{r+1}$, which corresponds to the $(r+p)$th column of $\bar{X}$, can be computed straightforwardly given the other rows of $L_{r+1}$ and $\bar{X}_{C_{r+1} C_{r+1}}$.  The same technique can be repeated to construct the Cholesky factor of $\bar{X}_{C_{r+2} C_{r+2}}$ from $L_{r+1}$ and so on. This technique computes $L_{r+1}$ in $O(p^2)$ operations (as opposed to $O(p^3)$ operations if $L_{r+1}$ is computed from scratch) and hence reduces the total time to compute $\ln\det\hat{X}$ to $O(np^2)$ operations, which is the same order as the complexity of Cholesky factorization.  Also note that, incidentally, $R^T$ is the Cholesky factor of $\bar{X}_{U_r U_r}$. This coincidence does not always occur in general case.

This Cholesky updating process is not limited to the banded matrix case. In general, given the Cholesky factor $L_r$ of $\bar{X}_{C_r C_r}$ (resp. $\bar{X}_{U_r U_r}$), the Cholesky factor of a symmetric permutation of $\bar{X}_{C_{r'} C_{r'}}$ (resp. $\bar{X}_{U_{r'} U_{r'}}$) can be constructed by removing the rows of $L_r$ corresponding to $C_r \setminus C_{r'}$ (resp. $U_r \setminus U_{r'}$), performing Givens rotation to transform its transpose into an upper triangular matrix, and then appending the rows corresponding to $C_{r'} \setminus C_r$ (resp. $U_{r'} \setminus U_r$). The resulting matrix may not be the Cholesky factor of $\bar{X}_{C_{r'} C_{r'}}$ (resp. $\bar{X}_{U_{r'} U_{r'}}$) but rather of some symmetric permutation of it because the rows corresponding to $C_{r'} \setminus C_r$ (resp. $U_{r'} \setminus U_r$) are always appended to the bottom.  Since permuting a matrix symmetrically does not affect its determinant, the resulting Cholesky factor can be used for determinant computation as is. 

Roughly speaking, the above technique is more efficient the smaller $|C_r \setminus C_{r'}|$ and $|C_{r'} \setminus C_r|$ are.  Larger $|C_r \setminus C_{r'}|$ \emph{usually} implies more Givens rotations while larger $|C_{r'} \setminus C_r|$ implies more computation of the entries of the appending rows.  However, if the rows to be removed are the bottom rows of $L_r$, no Givens rotations are required (as the resulting matrix remains lower triangular). Therefore, large $|C_r \setminus C_{r'}|$ does not imply many Givens rotations in this case.

Hence, we are able to compute $\ln\det\hat{X}$ ``optimally" (in the sense of within the same order of complexity as performing Cholesky factorization) in two special cases: block diagonal and banded matrices. However, there are cases where it seems $\ln\det\hat{X}$ cannot be computed ``optimally" using Cholesky updating scheme, for example, the planar graph case.  It is this reason that prevents us from using the dual algorithm described in Section \ref{section_dualnewton} to compute the primal projected Newton direction as it requires evaluation of $(\nabla^2\ln\det\hat{X})P$ in each iteration of the conjugate gradient and evaluating $(\nabla^2\ln\det\hat{X})P$ is generally at least as expensive as evaluating $\ln\det\hat{X}$.

\subsection{Computation of $\hat{X}^{-1}$ and $(\nabla^2\ln\det\hat{X})P$}
\label{section_xderiv}

As mentioned in Section \ref{section_primalnewton}, computing $N$ in step (ii) of DPD requires the knowledge of $\hat{X}^{-1} = \nabla\ln\det \hat{X}$. In addition, steps (iii), (iv), and (ix) also call for $\nabla\ln\det\hat{X}$ and $(\nabla^2\ln\det\hat{X})N$ in the formula for $\Delta S_1$ and in the steepest descent direction, respectively. From (\ref{logmaxdet}), the matrix $\hat{X}^{-1}$ is seen to be 
\begin{eqnarray}
\label{logmaxdetderiv}
\hat{X}^{-1} & = & \frac{d}{d \bar{X}}\ln\det \hat{X} \nonumber \\
&=& \sum_{r=1}^l \left(\frac{d}{d\bar{X}}\ln\det \bar{X}_{C_r C_r}\right) - \sum_{r=1}^{l-1} \left(\frac{d}{d\bar{X}}\ln\det \bar{X}_{U_r U_r}\right). 
\end{eqnarray}
Recall that $\bar{X}_{C_r C_r} (r=1,2,\ldots,l)$ and $\bar{X}_{U_r U_r} (r=1,2,\ldots,l-1)$ are completely dense.  For this reason, using automatic differentiation would not yield a more efficient first-derivative computing algorithm than simply computing $\frac{d}{d\bar{X}}\ln\det \bar{X}_{C_r C_r}$ and $\frac{d}{d\bar{X}}\ln\det \bar{X}_{U_r U_r}$ conventionally (by finding their inverses) and piecing them together according to (\ref{logmaxdetderiv}).  The same is true with the product of the second derivative and an arbitrary matrix
\[
(\frac{d^2}{d\bar{X}^2}\ln\det\hat{X})P' = \sum_{r=1}^l \left(\frac{d^2}{d\bar{X}^2}\ln\det \bar{X}_{C_r C_r}\right)P' - \sum_{r=1}^{l-1} \left(\frac{d^2}{d\bar{X}^2}\ln\det \bar{X}_{U_r U_r}\right)P'.
\]
Recall that $\left(\nabla^2\ln\det W \right)P' = -W^{-1}P'W^{-1}$.  Therefore, our algorithm computes $\frac{d}{d\bar{X}} \ln\det\hat{X}$ by the simple algorithm described above.  The computation of $(\frac{d^2}{d\bar{X}^2}\ln\det\hat{X})P'$ is also handled similarly: by computing the product of the second derivative of each dense submatrix and $P'$ conventionally 
and then piecing them together.

The Cholesky updating technique as described in Section \ref{section_maxdet} can also be applied to the computation of $\frac{d}{d\bar{X}} \ln\det\hat{X}$ and $(\frac{d^2}{d\bar{X}^2}\ln\det\hat{X})P'$. Both of these computations involves computing the Cholesky factor of each clique in order to compute its inverse.  Hence, the same technique can be applied to reduce the computation time required to find the Cholesky factors.

\section{Estimates of time and space complexities}

We give estimates of time and space complexities of our algorithm in this section. Let $O(\mathrm{Time}(L_F))$ and 
$O(\mathrm{Space}(L_F))$ denote the time and space complexity of Cholesky factorizing a matrix 
with sparsity pattern $F$, respectively. Steps (i) and (v) do not require
any computation. Step (ii) involves a conjugate gradient to solve (\ref{lag2}). Each iteration of the conjugate gradient requires one evaluation of the Hessian-vector product $(\nabla^2 q(\mathbf{0}))\mathbf{\lambda}$, 
which is $O(\mathrm{Time}(L_F))$ and $O(\mathrm{Space}(L_F))$. The conjugate gradient 
takes at most $n$ iterations to converge. Step (iii) requires one evaluation of 
$(\nabla^2\ln\det\hat{X})N$. Step (iv) requires one evaluation of $\nabla\ln\det\hat{X}$ and $(\nabla^2\ln\det\hat{X})N$ each. Step (vi) involves a conjugate gradient that takes one evaluation of the Hessian-vector product $(\nabla^2 h(\mathbf{y}))\mathbf{z}$ and therefore requires $O(\mathrm{Time}(L_F))$ and 
$O(\mathrm{Space}(L_F))$ per conjugate gradient
iteration. The conjugate gradient also takes at most $n$ iterations to converge. Step (vii) requires one evaluation of 
the entries in $F$ of $(\nabla^2\ln\det S)\tilde{N}$, which is $O(\mathrm{Time}(L_F))$ and $O(\mathrm{Space}(L_F))$. Step (viii) calls for
one evaluation of the entries in $F$ of $\nabla\ln\det S$ and $(\nabla^2\ln\det S)\tilde{N}$ each. Finally, step (ix)
requires one evaluation of $\ln\det\hat{X}$ per iteration of steepest descent.

A few steps in the algorithm, namely steps (iii), (iv), and (ix), require evaluation of either $\ln\det\hat{X}$ or one of its derivatives. Generally, evaluating $\ln\det\hat{X}$ or its derivative is more expensive than $O(\mathrm{Time}(L_F))$. Fortunately, the computational results in Section \ref{section_exp} show
that only a constant number of steepest descent iterations are needed to find good step size, and steps (iii) 
and (iv) only require at most two evaluations of such quantities.  Space complexity of an evaluation of $\ln\det\hat{X}$,
on the other hand, is still $O(\mathrm{Space}(L_F))$ as evaluating $\ln\det\hat{X}$ reduces to computing the Cholesky factor
of each of the maximal cliques without having to store the Cholesky factor of more than one clique at a time. 
For the case where $G=(V,E)$ is planar and $F$ is its chordal extension when the vertices $V$ are ordered in nested dissection ordering, computing $\ln\det\hat{X}$ (and, consequently, each iteration of the steepest descent method) requires $O(n^{5/2})$ operations and $O(n \log n)$ space. Notice that $O(n^{5/2})$ operations for these steps are acceptable as 
each of the conjugate gradient takes $O(\mathrm{Time}(L_F)) = O(n^{3/2})$ operations per iteration and 
at most $n$ iterations to converge, resulting also in $O(n^{5/2})$ operations for steps (ii) and (vi).

As the last remark for this section, we note that we suspect that this estimate of $O(n^{5/2})$ for 
planar case may not be tight. For the special case that $G$ is a grid graph, it is not hard to show that computing $\ln\det\hat{X}$ requires only $O(n^{3/2})$ operations, which is the same order
 as performing Cholesky factorization.

\section{Computational results}
\label{section_exp}

We implemented and tested our algorithm by using it to solve various instances of the problem of finding maximum cut (MAX-CUT). The procedure of using SDP to solve MAX-CUT is proposed by Goemans and Williamson in 1995 \cite{goemans}. Readers are referred to Goemans and Williamson's paper for the details of the procedure. In MAX-CUT, the input graph whose maximum cut is sought is exactly the aggregated sparsity pattern $E$ of the resulting SDP program.

Given the aggregated sparsity pattern $E$, we find its chordal extension by ordering the vertices of $E$ according to the symmetric minimum degree ordering \cite{markowitz}, perform symbolic Cholesky factorization on the reordered matrix, and use the resulting Cholesky factor as the chordal extension. The primal partial variable is initialized to the identity matrix. The dual variable is initialized to $C-I$ after $C$ has been reordered according to the minimum degree ordering. After the algorithm finds an iterate whose duality gap is less than $10^{-3}$, it continues for 3 additional iterations and then terminates. Each conjugate gradient runs until $\|r\|_2$ is less than $10^{-5}$ times the 2-norm of the constant term of the system that the algorithm is trying to solve. Finally, step (ix) of the algorithm is implemented using the method of steepest descent starting from four initial points $(X+\Delta X_1,S)$, $(X+\Delta X_2,S)$, $(X,S+\Delta S_1)$, and $(X,S+\Delta S_2)$ separately (Refer to chapter 6.5.2 of \cite{heath} for the explanation of the method of steepest descent). We do not perform any line searches in the steepest descent; we simply take the step size to be identically one and take the step as long as the new point decreases the potential. Line searches are ignored because, according to our testing, performing 
line searches does not generally improve computation time of the algorithm. 
The additional evaluation of $\ln\det\hat{X}$ in each step of the line searches 
appears to be too expensive compared to the extra decrease in potential resulted from them.

The test instances were generated by adding edges to the graph randomly until the chosen number of edges were met. The number of main iterations reported in column 5 of Table \ref{table_main} is the number of times the algorithm repeats step (i) to (x) before it finds an optimal solution is found. 

\begin{table}
\begin{center}
\begin{tabular}{|r|r|r|r|r|r|r|r|}
\hline
$n$ & $m$ & Trials & Time (s)& Num Iter & $\Delta X_1$ CG & $\Delta S_2$ CG & Pot Min\\ 
\hline \hline
5   & 7   & 100 & 1.30   & 13.8 & 3.8  & 3.8  & 2.0 \\
10  & 16  & 100 & 5.43   & 15.8 & 8.8  & 9.2  & 2.2 \\
20  & 40  & 50 & 25.12  & 17.9 & 14.9  & 15.4 & 2.4 \\
50  & 75  & 10 & 124.33  & 21.9 & 24.9  & 26.1 & 3.6 \\
100 & 180 & 5  & 894.44 & 25.2 & 35.3 & 37.0 & 4.9 \\
\hline
\end{tabular}
\end{center}
\caption{Summary of results of the algorithm on random instances. From left to right, the columns are the number of vertices, the number of edges, the number of trials run, the average CPU time in seconds, the number of main iterations, the average number of conjugate gradient iterations required to compute the search direction $\Delta X_1$ for one point, the average number of conjugate gradient iterations required to compute $\Delta S_2$ for one point, and the average iterations to minimize potential along the four directions for one starting point (step (ix) in the algorithm).} 
\label{table_main}
\end{table}

The results in Table \ref{table_main} show that the conjugate gradient to compute $N$, the matrix necessary for computation of $\Delta X_1$, requires more number of iterations than to compute $\tilde{N}$, the matrix necessary for computation of $\Delta S_2$. In addition, the potential minimization by the method of steepest descent only takes a few iterations for each initial point and therefore does not steal away too much valuable time that could be spent on other computations.

The two directions $\Delta X_2$ and $\Delta S_1$ are not common in literature although the other two directions $\Delta X_1$ and $\Delta S_2$, which are the projected Newton directions, appear in many other SDP algorithms. We performed an experiment to test whether using all 4 directions are more beneficial than using only 2 more common directions. We tested the two versions of our algorithms on the random instances generated in the same manners as the ones in the previous experiment. The results are shown in table \ref{table_4vs2}.

\begin{table}
\begin{center}
\begin{tabular}{|r|r|r|r|r|}
\hline
$n$ & $m$ & Num Trials & \multicolumn{2}{|c|}{Time (s)} \\ 
\cline{4-5}
& & & 4 directions & 2 directions \\ 
\hline \hline
5 & 7 & 100 & 1.30  & 1.23\\   
10 & 16 & 100 &5.43 & 7.58\\
20 & 40 & 50 &25.12 & 37.40\\
50 & 75 & 2 &124.33 & 191.46\\
\hline
\end{tabular}
\end{center}
\caption{Comparison of the average CPU time (in seconds) the two algorithms required to solve the random instances.} 
\label{table_4vs2}
\end{table}

Table \ref{table_4vs2} shows that using all four directions make the algorithm find the optimal solution in shorter time in all test cases. The reason toward this result is that all of the quantities involved in computation of the two uncommon directions are also required to compute the other two projected Newton directions.  Therefore, computation of the additional two unusual directions is relatively cheap compared to the reduction in potential they induce.

\begin{figure} 
 \centering 
 \includegraphics[width=0.9\textwidth]{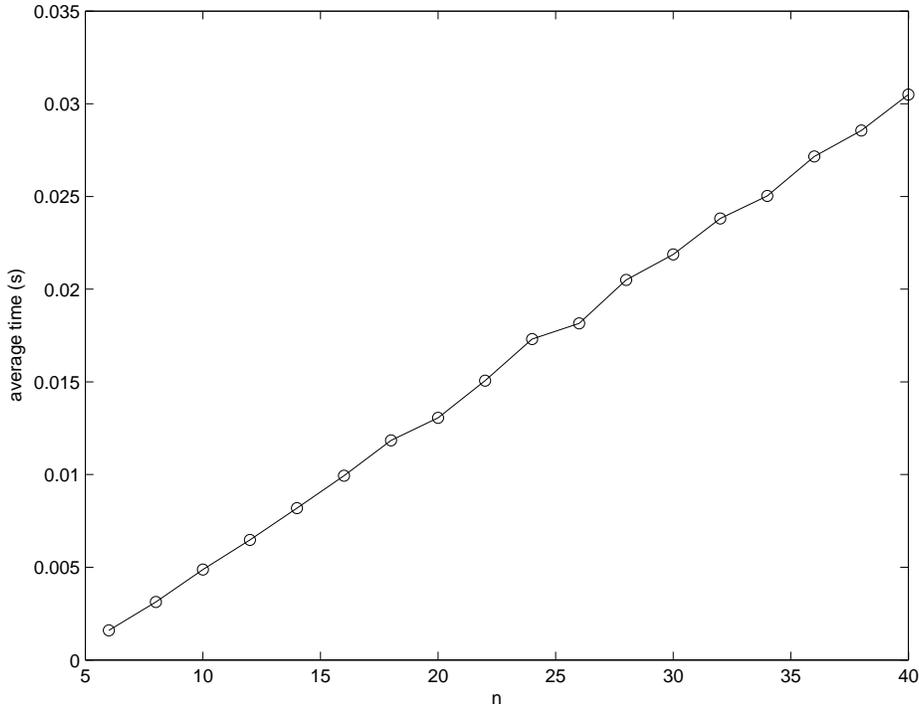} 
 \caption{Average CPU time to compute $\ln\det\hat{X}$ in the case that $\bar{X}$ is a banded matrix. Bandwidth is fixed to 3 while varying the number of vertices.} 
 \label{figure_band_vary_n} 
\end{figure}

\begin{figure} 
 \centering 
 \includegraphics[width=0.9\textwidth]{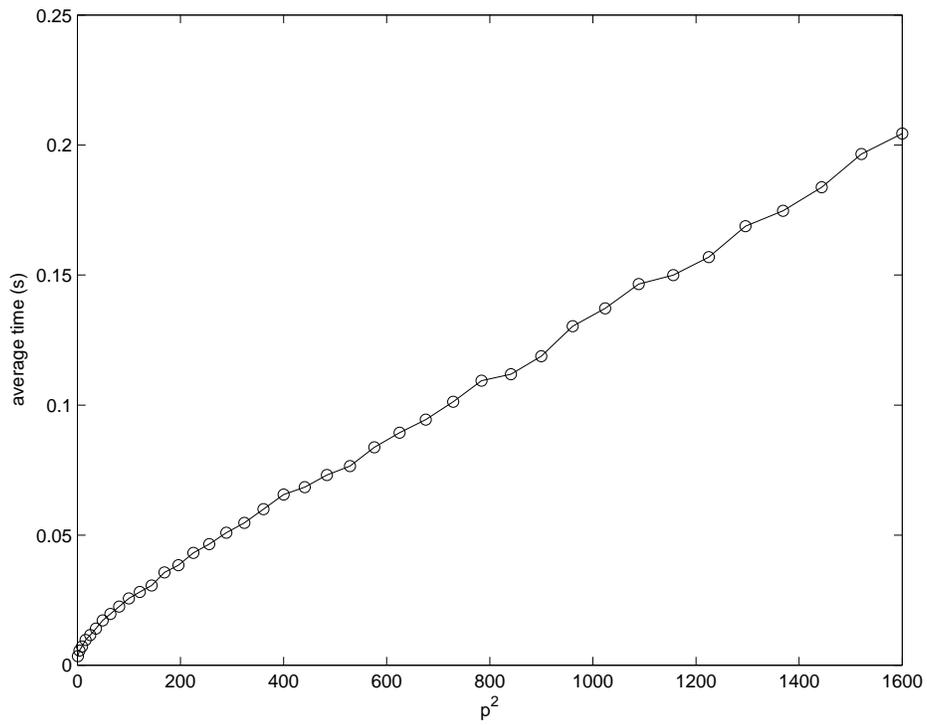} 
 \caption{Average CPU time to compute $\ln\det\hat{X}$ in the case that $\bar{X}$ is a banded matrix. The quantity $n-p$ is fixed to 10 while varying the bandwidth.} 
 \label{figure_band_vary_p} 
\end{figure}

Finally, we tested the computation of $\ln\det\hat{X}$ with Cholesky updating scheme on banded matrices to verify its $O(np^2)$, or more precisely, $O((n-p)p^2)$ complexity. We began by fixing the bandwidth $p$ to be 3 and varying $n$ from 6 to 40, repeated 500 times for each $n$. Figure \ref{figure_band_vary_n} shows the plot of the average CPU time to compute $\ln\det\hat{X}$ for banded matrices with bandwidth 3 of various size against the number of vertices, and it confirms the linear dependency on $n$ of the complexity. Next, we fixed $n-p$ to 10 and varying $p$ from 1 to 40, repeated 50 times each. The plot of the average CPU time against the square of the bandwidth for this experiment is shown in figure \ref{figure_band_vary_p}. The plot agrees that the complexity of $\ln\det\hat{X}$ is proportional to $p^2$ in the banded case. 


\section{Concluding Remarks}

We showed an implementation of a SDP solver that exploits sparsity in the data matrices for both primal and dual variables. Our algorithm is based on the primal-dual potential reduction method of Nesterov and Nemirovskii and uses partial primal matrix variable as proposed by Fukuda et al. Two of the search directions are projected Newton directions that can be found by solving linear systems involving the gradient and the Hessian of the logarithm of the determinant of a matrix with respect to a vector. We observed that the idea from reverse mode of automatic differentiation can be applied to compute the mentioned gradient and the product of the Hessian and an arbitrary vector efficiently, which is in the same order as computing determinant of a sparse matrix. Using this observation, we solve the linear system for the search directions by conjugate gradient, which requires one evaluation of the product of the Hessian and a vector in each iteration in exchange for not having to factorize the Hessian matrix. For the primal case, we propose a way to compute one of the primal search directions without requiring the determinant of the positive definite completion in each iteration of the conjugate gradient because
the determinant of such completion is generally more expensive than performing Cholesky factorization. This determinant is still required in the potential minimization to find step sizes as well as in the course of computing one other search direction. we described a technique to reduce the complexity of computing the logarithm of the determinant of a positive definite matrix completion by reusing the Cholesky factors when there are many overlaps of maximal cliques. This technique reduces the complexity to that of performing Cholesky factorization in the banded matrix case but still cannot achieve the same complexity as the Cholesky factorization in general. Fortunately, only a few number of evaluations of the determinant of such completion is required per one SDP iterate. The other two non-Newton directions can be computed efficiently since they require the same quantities that are already computed in the process of finding the former projected Newton directions. We then tested our algorithm on random instances of the MAX-CUT problems. From the results, the conjugate gradients do not require too many iterations to converge for the algorithm to be impractical.

There are questions unanswered in this paper that can help improve the algorithm described here. For example, can we compute the logarithm of the determinant of a positive definite matrix completion more efficiently, perhaps by another derivation different from (\ref{cliquefactor})? The other issue is regarding the stability. How can we incorporate preconditioning to mitigate the ill-conditioning of the linear system problems? Regardless, our algorithm should prove efficient in the applications where the data matrices are sparse.

\end{document}